# Optimal Price and Quantity Determination of Retailer Electric Contract and maximizing social welfare in Retail Electrical Power Markets with DG

Barati Masoud, Nikkhah Mojdehi Mohammad, Kazemi Ahad

*Center of Excellent for Power System Automation and Operation Department of Electrical Engineering, Iran University of Science & Technology*

*Abstract*—Retail electrical power marketers, also known as retailers, typically set up contracts with suppliers to secure electricity at fixed prices on the one hand and with end users to meet their load requirements at agreed rates on the other hand. Also high penetration of distributed generation (DG) resource is increasingly observed worldwide. Considering the viewpoint of a retailer, this paper analyzes the problem of setting up contracts on both the supplier and end-user sides to maximize profits while maintaining a minimum operational cost of distribution system. The proposed two-level optimization models can minimize the cost of distribution system with Distributed Generation and maximize the profits of retailers. The numerical simulations are carried out based on an IEEE 33-bus test distribution network. Test results are included to show the performance of the proposed method.

*Index Terms*—Retailer profit, electricity market, wholesale market, Distributed Generation (DG), social welfare, Distribution system, DisCo, distribution system operator (DSO).

## I. INTRODUCTION

THE U.S. electric power market has lately undergone significant changes due to deregulation and restructuring of the industry. One result of this new marketplace is the emergence of third-party entities known as retailers or marketers. These entities purchase power from suppliers and sell it again to end-user customers. These retailers typically set up bilateral contracts with suppliers to procure a substantial portion of the required power to meet the demands of their customers. The remaining portion of the load is bought in the balancing energy market and is settled at the marginal clearing price or the spot price. Retailers set up supply contracts with suppliers. These contracts range from a simple fixed-price, fixed-quantity contract for a fixed time period—where the retailer takes all the risk-or to a more collaborative arrangement where the supplier and the retailer work on a profit. Retailers are also turning to other market instruments to

Manage their portfolio and their risks. They do this by purchasing and selling options, hedging portions of their needs especially during the volatile months, and using cross-commodity hedges (natural gas as a hedge for electricity).

Distributed generation (DG) can be defined as the integrated use of small generation units directly connected to a distribution system or inside the facilities of a customer [1]. The potential development of DG is sustained in the following factors: increasing power quality requirements, avoiding or shifting investment in transmission lines and/or transformers, minimizing ohmic losses, protecting the environment, and maintaining high energy prices at retail level [1]–[4].

Traditionally, a distribution company (DisCo) purchases energy from the wholesale market, at a high voltage level, and then transfers this energy to final customers. Nevertheless, the restructuring process of the energy sector has stimulated the introduction of new agents such as retailers and products, and the unbundling of traditional DisCo into technical and commercial tasks, including the provision of ancillary services [2].

While there has recently been extensive research on market modeling at the wholesale level, the retail side and the reaction with DG has received relatively less attention.

Several researches have been carried out to find a way to the retail competition and reaction with distributed generation problems on the retail market. In Reference [5] several strategies have been analyzed for determining forward loads for electrical power retailers based on empirically fit probability distributions for end-user loads and spot market prices. In this paper it is also assumed that in a distribution network the retailer contracts with classes of the final consumers. In reference [6] there have been addressed the contract design problem a retailer may encounter, both at the supply and end-user levels. It provides a stochastic programming methodology that allows the retailer to make informed contractual decisions, particularly in respect to contract prices and quantities. In accordance with [6] it should be emphasized that an appropriate modeling of the risk associated with buying from the spot market is required. In reference [7] a model has been defined that can be applied to an energy service provider (which supplies electric energy for eligible industrial customers) in the Spanish electricity market. In accordance with [7] this model results in a quadratic non-linear optimization problem, where the objective function to be maximized is the economic profit obtained by the supplier, and it calculates the optimal electricity selling prices to be applied to the customers. The novel market integration mechanism has been presented in [8] which demonstrate that the design of a market interface and the consideration of DG units as equivalent power producers behave in accordance with the wholesale market methodology. The main concepts of the proposal are presented theoretically and are illustrated.

The authors are with the Center of Excellent for Power System Automation and Operation Department of Electrical Engineering, Iran University of Science & Technology, Narmak 16846-13114 Tehran, IRAN (e-mail: m.barati@moe.org.com; nikkhah627@gmail.com; kazemi@iust.ac.ir).

In accordance with [8] the proposed scheme shows an adequate trade-off between technical accuracy (marginal cost theory and ohmic losses) and practicability (available network information and information management).

A Multiperiod energy acquisition model for a DisCo in the competitive day-ahead electricity market with two new resources, distributed generation and interruptible load have been proposed in [9]. Due to the interactions among Discos' strategies, the energy purchasing scheme is modeled as a bilevel optimization problem in order to maximize every Disco's profits. In accordance with [9] to solve the complementary problems formulated from the KKT conditions, a nonlinear complementary method has been employed to solve the proposed model.

The goal of this paper is to find the sale price of retailers in the way that not only increases the profit of retailers but also develops social welfare. By social welfare we mean minimizing the cost of putting the distribution system into operation. This goal is implemented in a distributional network which in itself contains some distributed generation resources in a way that affects the mechanism of sale business. The costs will include the production cost of distributed generation units and the cost of purchasing from the wholesale market. To answer the above problem, an algorithm has been proposed which is implemented in two phases.

The rest of this paper is organized as follows. In Section II, we describe modeling considerations and provide mathematical formulations for the simulation. Section III describes computations and analyses presented in Section II; Section IV describes our numerical findings and is followed by the conclusions section.

## II. RETAIL ELECTRIC POWER MARKET MODEL

The model which is proposed in this section, considers two principles.
- The first one concerns minimizing the cost of distribution system operation, which is a problem Disco always faces, which means how one should use available resources to minimize the cost of distribution system operation. In the proposed model, it is assumed that available resources for Disco are distributed generation and purchased energy from retailers.
- The second principle is about maximizing the retailer's profit. It must be considered that the retailer's profit be maximized. Retailers, on one hand purchase power from wholesale market and distributed generation resources with spot prices considering the marginal cost of distributed generation units respectively and on the other hand sell power to their final consumers with their expected sale price.

The combination of these two pre-mentioned points formulates an important problem in the way of achieving the expected results. The proposed method in this paper is based on a loop.

Meanwhile the goal of this paper is to find the sale price of retailers in the way that not only increases the retailer's profit but also develops social welfare. By social welfare we mean minimizing the cost of distribution system operation. This goal is implemented in a distribution network which in itself contains some distributed generation resources in a way that affects the mechanism of sale business. The costs will include the production cost of distributed generation units and the cost of purchasing from the wholesale market.

So the profit of any retailers will be mentioned this way:

$$\text{Profit} = \text{income} - \text{payment} \quad (1)$$

In (1), income includes the input of each retailer based on the existing load in its class and the tariff on which the contract is based. As it was mentioned in section I in this paper it is also assumed that in a distribution network the retailer contract with classes of the final consumers. These classes are usually based on their geographical situations. Payment will include the payout of each retailer to the power supplier with whom the named retailer contracts to provide the load of existing consumers in its own class. The resources from which the retailer in the distribution network can provide the power are the wholesale market and distributed generation resources. It is obvious that if there is no DG in the distribution network, the retailer will provide all the required power from the wholesale market with spot prices. But if there is any distributed generation units in the distribution network, then we will face a critical question and that is "what amount of required power would a retailer provide from the wholesale market and what amount will be provided from distributed generation units?" As a matter of principle the combination of the retailer's purchase from the wholesale market and the distributed generation units would be in a way to increase the gaining profit by the retailer. But would this be advantageous for the final consumer to answer all these questions, in the following an algorithm has been proposed to reach the maximum social welfare.

## III. MATHEMATICAL FORMULATION

As it was mentioned earlier, the presence of distributed generation resources in the distribution network will bring about some ambiguities concerning the payout of the retailer and the purchasing strategy carried out by them. To remove these ambiguities an algorithm is proposed in this paper that will be accomplished in two phases. At the first phase, it is trying to answer the question that by what strategy should a retailer purchase the required power with spot price from the wholesale market and the distributed generation units that would both maximize his profit and minimize the cost of operating the whole distribution system and therefore increase the social welfare of the system. In this way at the first phase of this proposed algorithm an optimal power flow will be administered from the point of view of the distribution system operator (DSO), in a way that the cost of power production in the distribution system must be minimized. By using this method the production portion of each distributed generation unit and the power injection portion from the wholesale market to the distribution network classes will be determined. The objective function of this optimal power flow will be as

follows:

$$Min \sum_{c=1}^{n} \frac{Load_c}{\sum_{c=1}^{n} Load_c} \cdot P_W \cdot Price_c + \sum_{i=1}^{m} C_{DGi}(P_i) \quad (2)$$

Subject to:

$$P_W + \sum_{i=1}^{m} C_{DGi}(P_i) = P_{loss} + \sum_{c=1}^{n} Load_c \quad (3)$$

$$V_i^{min} \leq V_i \leq V_i^{max} \quad (4)$$

$$P_i^{min} \leq P_i \leq P_i^{max} \quad (5)$$

Where $n$ denotes the number of classes; $Load_c$ denotes real power demand of class $c$; $P_W$ denotes the Injection power from the wholesale market; $Price_c$ denotes the price of electrical energy class $c$; $P_{loss}$ denotes real power loss in distribution network; $m$ denotes the number of DG units; $C_{DGi}(P_i) = a_i P_i^2 + b_i P_i + c_i$ denotes the cost characteristic of DG at node $i$.

The first term in the objective function (2) represents the costs of energy purchasing from wholesale market, and the second term represents the generation costs for distributed generation units. Constraints (3), (4) and (5) represent the load balance of the whole distribution network constraint. Generating limit are specified as upper and lower limits for the real power outputs of DG unit and the voltage constraint on every node of the distribution network respectively.

For the beginning of implementing this algorithm at the first phase, we do not have any access to the price of classes; therefore we need to make some initial values to determine the price of the classes. At the first iteration it is assumed that the price of the classes equals the spot price of the wholesale market, so a primary value for the generation of distributed generation units and the injection power from the wholesale market will be achieved.

Using these initial values, the algorithm will enter its second phase. At this phase, the main goal is to maximize the profit of retailers in a way that according to (2) the income of purchasing electricity to classes be maximized but on the other hand the costs related to purchasing electricity from the wholesale market with spot price and purchasing from scattered resources with limited prices of these production units be minimized. By fulfilling this part of algorithm the price of electricity will be measured in each class.

The objective function in this part would be as follows:

$$Max \sum_{c=1}^{n} Load_c \cdot Price_c - \sum_{c=1}^{n}\sum_{i=1}^{m} P_{c\_DGi} \cdot MC_{DGi} - Price_{spot} \cdot \sum_{c=1}^{n} P_{c\_W} \quad (6)$$

Subject to:

$$Load_c = Load_c^N \cdot \left[1 + \frac{\beta_c(Price_c - Price_c^N)}{Price_c^N}\right] \quad (7)$$

$$Load_c = P_{c-W} + \sum_{i=1}^{m} P_{c-DGi} \quad (8)$$

Where $P_{c\_DGi}$ denotes the Power purchased by a retailer of class $c$ from the distributed generation unit $i$ ($DG_i$); $MC_{DGi}$ denotes the marginal cost of distributed generation unit $i$; $Price_{spot}$ denotes the Spot price of the wholesale market; $P_{c\_W}$ denotes the power purchased by a retailer of class $c$ from the wholesale market; $Load_c^N$ denotes the nominal load of class $c$; $\beta_c$ denotes the Load elasticity class $c$; $Price_c^N$ denotes the nominal sale price of class $c$; $Price_c$ denotes the sale of class $c$.

The first term in the objective function (6) represents the income of energy selling of each retailer correspondent to their classes, the second term represents the purchasing energy from distributed generation units and the third term represents the purchasing energy from the wholesale market. Constraints (7) and (8) represent the end user demand function and the load balance of each class respectively.

In [10, 11] one can observe one of the existing and most useful models of the function of price response of loads of network. In (7) the load function has come in linear form based on nominal variables. In this equation $Price_c^N$ is the selling price at the nominal electrical power consumed. In order to calculate the nominal selling price $Price_c^N$ is applied to the customers. Also it is assumed that (in the situation of nominal demand $Load_c$) the retailer sells electricity at a price which is 10% higher than the whole cost, that is, it assigns a profit on income of 10%.

In this equation it is assumed that the response of each consumer to the price variation, characterized by the elasticity $\beta$, has to be empirically assigned. The price elasticity of demand $\beta$ measures the degree of response of the demand to the changes of price in the market, being defined as the ratio of the variation of the demand quantity of some commodity or service, in percent, to the variation of its price $Price$, also in percent. It can be expressed as:

$$\beta = \frac{\Delta Load/Load}{\Delta Price/Price} \quad (9)$$

The demand is 'elastic' if $\beta$ is higher than 1, and 'inelastic' if $\beta$ is lower than 1. Considering various reference studies performed to characterize the end user response to several price strategies [12-19], elasticity values ranging from -0.01 to -0.25 have been chosen. Corresponding to Fig.1 elasticity value is determined by a negative value and equals to -0.2.

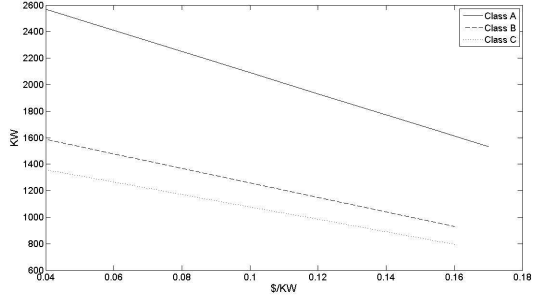

Fig.1. Demand function

After implementing of optimization problem in phase two, the sale price in each class will be achieved. Then these results will be applied to optimal power flow in phase one as input. After optimal power flow in phase one, the amount of economical generation from distributed generation units and amount of injection power from the wholesale market are determined. Then these economical power generation and injection from optimal power flow in phase one are applied to optimization problem of phase two as inputs in order to calculate sell prices. Until converging, the sale price of each class must be iterated. By using proposed algorithm, it can be seen that without load elasticity, retailer's sale price will increase sharply. On the other hand, load elasticity regulates sale prices and limits retailer's profit.

## IV. SIMULATIONS AND NUMERICAL RESULTS

An IEEE 33-bus distribution system [20] is used to illustrate the proposed model and solution algorithm. The system includes four distributed generation unit at nodes 4, 7, 25 and 30, three retailers; subsequently there will be three classes of retailing in the mentioned system. Table I shows the cost characteristics of varieties of DG unit technologies considered in this work. In this system the number of each load will be determined by the number of the node attached to it. Table II shows the load classification based on the classes.

TABLE I
DISTRIBUTED GENERATION DATA
WITH VARIETIES TECHNOLOGY [21]

| DG Technology | a | b | c | $P_{max}$ | $P_{min}$ |
|---|---|---|---|---|---|
| Fuel cell- CHP | 0.0001 | 0.0315 | 1.0749 | 400 | 0 |
| Gas ICE-CHP | 0.0001 | 0.0374 | 0.4777 | 400 | 0 |
| Gas ICE-power only | 0.0001 | 0.0777 | 0.3483 | 400 | 0 |
| Microturbine-CHP | 0.0001 | 0.0421 | 0.5553 | 400 | 0 |
| Microturbine-power only | 0.0001 | 0.0841 | 0.4603 | 400 | 0 |

Note: ICE: Internal Combustion Engine, CHP: Combined Heat and Power

TABLE II
LOAD CLASSIFICATION BASED ON RETAILER CLASS

| Retailer's Class | Loads |
|---|---|
| A | 2, 3, 4, 5, 6, 19, 20, 21, 22, 23, 24, 25 |
| B | 7, 8, 9, 10, 11, 12, 13, 14, 15, 16, 17, 18 |
| C | 26, 27, 28, 29, 30, 31, 32, 33 |

The assumed time interval for simulations is 24 hours a day. Fig. 2 and Fig. 3 show the load variations and spot prices of the wholesale market during a day respectively.

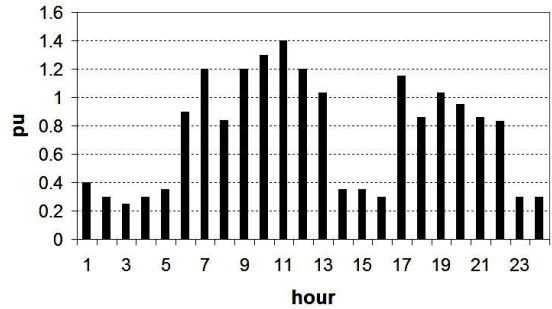

Fig. 2 Daily Load Variations

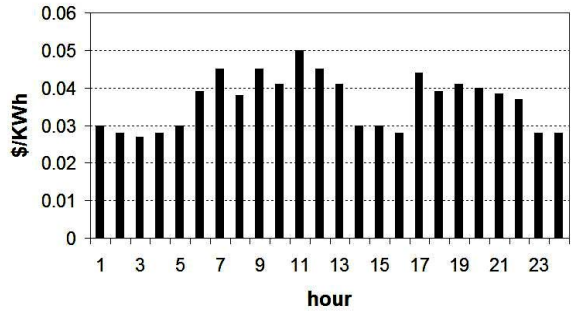

Fig. 3 spot prices of the wholesale market during a day

From Fig.4 to Fig.6, it can be seen that the retailer of each class tends to supply the demand of its class from the wholesale market principally. Also it can be concluded that the retailer's strategy is approximately constant with applying different distributed generation technologies and because of high consumed power in each class, retailers prefer to use this units in off peak hours rather than peak hours.

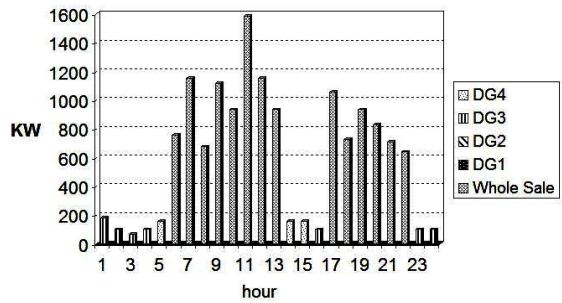

Fig. 3 purchasing strategy of retailer of class A, with Gas ICE-Power only technology

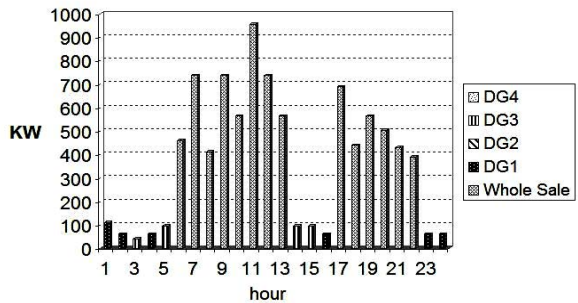

Fig. 4 purchasing strategy of retailer of class B, with Gas ICE- Power only technology

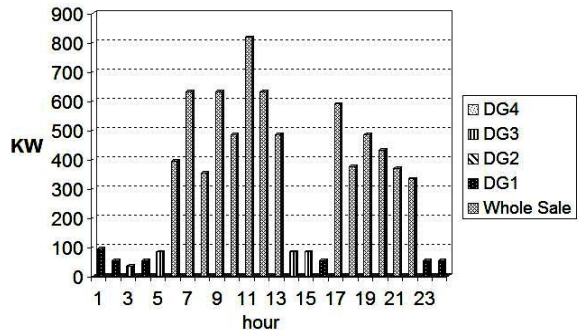

Fig. 5 purchasing strategy of retailer of class C, with Gas ICE- Power only technology

Fig.7 shows that the retailer's profit without using distributed generation resources is increased and out of evaluated technologies in this paper, Gas ICE-Power only and Fuel Cell CHP have the highest and the lowest profits for retailers respectively.

Also from Fig.7, it can be found that the retailer of class A has higher profit because of higher demand of its class than the others. More loads for the retailer dose not lead to more profit and the retailer's profit is affected by load forecasting. More loads are equal to more load forecasting error which can decrease the retailer's profit significantly.

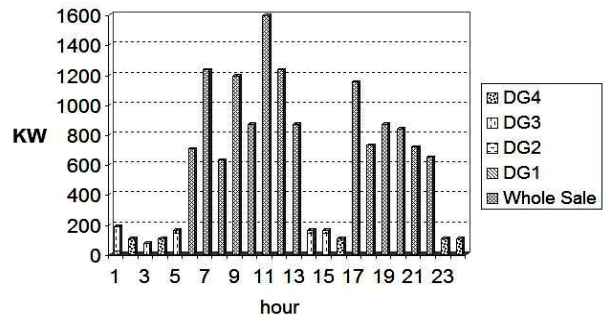

Fig. 6 purchasing strategy of retailer of class A, with Fuel Cell-CHP technology

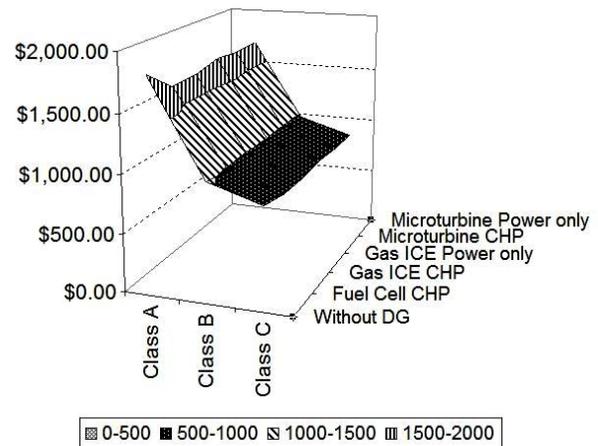

Fig.7 total profit of each retailer with several DG technologies

This issue has not been modeled in the proposed algorithm. Simulation results are come from this fact that the retailer's cost is composed of two parts. The first one is the cost of purchasing power from the wholesale market which is actually the purchased power multiplied by the price of the wholesale market which is constant during an hour and therefore the retailer's strategy relevant to this price will not be complicated. The second part of the retailer's cost is the cost of purchasing power from distributed generation resources which includes purchased power from these units multiplied by their marginal cost. It is worth to mention that the purchasing price (i.e. the marginal cost) from these units depends on its active generation power and is a linear curve with a positive slope. On the other hand, increase in power purchased from these units leads to higher price of these units therefore the cost of purchasing power from these units is increased. In comparison it should be mentioned that increase in purchasing power from the wholesale market leads to a lower purchasing power than the power supplied by distributed generation resources mostly because of constant prices.

An important point from these results is the retailer's income which includes the price of each class multiplied by the load of its class. As illustrated in Fig. 8 the retailer's price is constant during a period and dose not respond to hourly variations in wholesale prices and it results to uneconomical signal [21], So retailers look for suitable profit margin with lower cost and higher income.

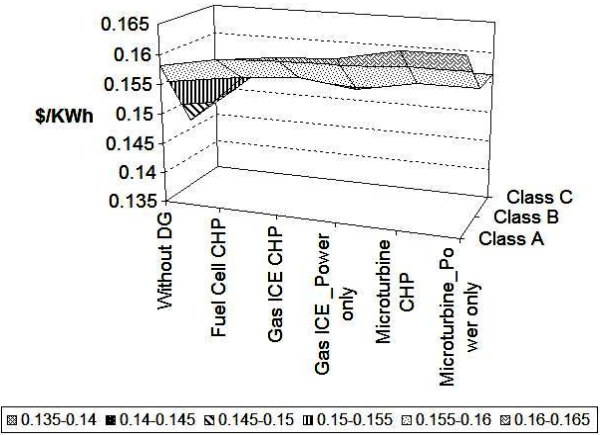

Fig.8. Price of each retailer with several DG technologies

## V. CONCLUSIONS

The following conclusions can be drawn from this paper:
- Retailers can experience the most profit without using distributed generation resources but their profits associated with distributed generation units are maximized with Gas ICE-Power only technology.
- In this paper, determination of optimal purchasing strategy and the retailer's profit maximization is implemented with social welfare maximization.
- Without applying subsidiary budget from governments, distributed generations can not be supported by retailers.
- In general terms, the more elastic the behavior of customer demand, the lower the profit obtained by retailers, as customers have an increased response capacity to price changes. On the other hand, most of the inelastic customers will tolerate price increases without reducing considerable consumption, and thus generate more revenues and higher profits.

## VI. REFERENCES

[1] T. Ackermann, G. Andersson, and L. Söder, "Electricity market regulations and their impact on distributed generation," in *Proc. Conf, Electric Utility Deregulation and Restructuring and Power Technologies 2000*, London, U.K., Apr. 4–7, 2000, pp. 608–613.
[2] "Organization for economic co-operation and development," in *Distributed Generation in Liberalized Electric Markets*: International Energy Agency, 2002.
[3] D. Kincaid. (1999) the Role of Distributed Generation in Competitive Energy Market. Gas Research Institute, Distributed Generation Forum. [Online]. Avalable: http://www.distributed-generation.com/library.htm
[4] A. P. Sakis Meliopoulos, "Distributed energy source: Needs for analysis and design tools," in *Proc. IEEE Vancouver Summer Meeting*, Vancouver, BC, Canada, 2001.
[5] S. A. Gabriel, M. F. Genc and S. Balakrishnan, "A Simulation Approach to Balancing Annual Risk and Reward in Retail Electrical Power Market," *IEEE Trans. Power system*, vol. 17, pp. 1050-1057, Nov. 2002.
[6] S. A. Gabriel, A. J. Conejo, M. A. Plazas and S. Balakrishnan, "Optimal Price and Quantity Determination for Retail Electric Power Contracts," *IEEE Trans. Power system*, vol. 21, pp. 180-187, Feb. 2006.
[7] J. M. Yusta, I. J. Ramirez-Rosado, J. A. Dominguez-Navarro and J. M. Perez-Vidal, "Optimal Electricity Price Calculation model for retailers in a deregulated market," *Elsevier. Electrical Power and Energy systems*, pp. 437-447, March. 2005.
[8] G. A. Jimenez-Estevez, R. Palma-Behnake ,R. Torres-Avila and L. S. Vargas, "A Competitive Market Integration Model for Distributed Generation," *IEEE Trans. Power system*, vol. 22, pp. 2161-2169, Nov. 2007.
[9] H. Li, Y. Li and Z. Li, "A Multiperiod Energy Acquisition Model for a Distribution Company with Distributed Generation and Interruptible Load," *IEEE Trans. Power system*, vol. 22, pp. 588-596, May. 2007.
[10] D. S. Kirschen, "Demand-side view of electricity markets," *IEEE Trans Power System,* vol. 18, pp. 520–527, 2003.
[11] R. N. Boisvert, P. A. Cappers and B. Neenan, "The benefits of customer participation in wholesale electricity markets," *Journal of Electra* vol. 15, pp. 41–51, 2002.
[12] P. Joskow and J. Tirole, "Retail electricity competition," CSEM working paper., University of California Energy Institute, Tech. Rep. WP-130, 2004.
[13] D. Royal (1164/2001). *Spanish Regulatory Legislation* [Online]. Available: http://www.cne.es
[14] F. C. Schweppe, M. C. Caramanis and R. D. Tabors. "Evaluation of spot price based rates," *IEEE Trans Power Appar System* vol. 7, pp. 1644–1655, 1985.
[15] F. Mochon, "Economı́a teo´rica y polı́tica," Madrid, McGraw-Hill; 1993.
[16] A. Tishler, "The industrial and commercial demand for electricity under time-of-use pricing," *Journal of Economics*, vol 23, pp. 369-384, 1983.
[17] D. J. Aigner and J. G. Hirschberg, "Commercial/industrial customer response to time-of-use electricity prices," some experimental results. *Rand J Econ*, vol 16, pp. 341-355, 1985.
[18] J. Zarnikau, "Customer responsiveness to real-time pricing of electricity," *Journal of Energy*, vol 4, pp. 99–116, 1990.
[19] J. C. Mak and B. R. Chapman, "A survey of current real-time pricing programs," *Journal of Electr*, vol. 7, pp. 76-77, 1993.
[20] J-H. Choi, J-C Kim, and S-I. Moon, "Integration Operation of Dispersed Generations to Automated Distribution Network for Network Reconfiguration," in *Proc. Conf, Power Technologies Conference 2003*, Bologna, Italy, June. 23-26, 2003.
[21] D. Holtz_Eakin, "Prospects for Distributed Electricity Generation," The Congress of the United States, Sep. 2003. [Online]. Available: http://www.cbo.gov

## VII. BIOGRAPHIES

**Masoud Barati** received the M.S degree in power system from Electrical Engineering Department of Iran University of Science and Technology (IUST) in 2005. Presently, he is an expert of Iran Electricity Market Regulatory Board (IEMR) in ministry of energy, Iran. His research interests are in market operation in electric power systems, Application of FACTS Devices to restructuring power system, Congestion management, power system dynamics and pss pricing, distributed generation, and Renewable energy.

**Mohammad Nikkhah Mojdehi** received the B.S degree in power system from Electrical Engineering Department of University of Guilan in 2006. Presently; he is student of Master Degree in Electrical Engineering Department of Iran University of Science and Technology (IUST). His research interests are in the application of heuristics method in power system, power system optimizations, and power system restructuring.

**Ahad Kazemi** received the M.S degree in electrical engineering from Oklahoma State University, U.S.A., in 1979. Currently, He is an Associated Professor at Electrical Engineering Department of Iran University of Science and Technology, Tehran, Iran. He is the chairman of department power system Iran University of Science and Technology, Tehran, Iran. His research interests are in the application of FACTS Devices in power system control, reactive power planning, and power system restructuring.